\documentclass[11pt]{article}

\usepackage{amsmath}
\usepackage{amssymb}
\usepackage{eucal}

\begin{document}

\baselineskip 16pt

\title{Finite groups with all $3$-maximal subgroups
$K$-$\mathfrak U$-subnormal}
\author{Xiaolan Yi \thanks{Research of the first author is supported by
a NNSF grant of China (Grant \# 11101369) and  the Science Foundation of Zhejiang
Sci--Tech University under grant 1013843-Y.}\\
{\small Department of Mathematics, Zhejiang Sci-Tech University,}\\
{\small Hangzhou 310018, P.R.China}\\
{\small E-mail: yixiaolan2005@126.com}\\ \\
{ Viktoria A. Kovaleva}\\
{\small Department of Mathematics,  Francisk Skorina Gomel State University,}\\
{\small Gomel 246019, Belarus}\\
{\small E-mail: vika.kovalyova@rambler.ru}}

\date{}
\maketitle

\begin{abstract} 
A classification of finite groups in which 
every $3$-maximal subgroup is $K$-$\mathfrak U$-subnormal is given.
\end{abstract}

 \let\thefootnoteorig\thefootnote
\renewcommand{\thefootnote}{\empty}
\let\thefootnote\thefootnoteorig

\section{Introduction}

Throughout this paper, all groups are finite and $G$ always denotes
a finite group. 
We use $\mathfrak U$
to denote the class of all supersoluble groups;
$G^{\mathfrak U}$ denotes the intersection of all normal subgroups 
$N$ of $G$ with  $G/N \in \mathfrak U$.
The symbol 
$\pi (G)$  denotes  the  set  of  prime  divisors  of 
the order of $G.$

A subgroup $H$ of  $G$ is called a {\sl $2$-maximal}  
({\sl second maximal}) subgroup of $G$ whenever $H$ is 
a maximal subgroup of some maximal subgroup $M$  of $G.$ 
Similarly we can define {\sl $3$-maximal subgroups}, 
and so on.

One of the  interesting  and  substantial  direction  in  finite  group  theory  
consists  in  studying  the  relations 
between the structure of the group and its 
$n$-maximal subgroups.
The earliest publications  
in  this  direction  are  the  articles  of  L. R\'{e}dei \cite{redei} 
and B. Huppert \cite{huppert}.
L. R\'{e}dei  described the nonsoluble groups with abelian second maximal 
subgroups.
B. Huppert  established  the  supersolubility  of    
$G$ whose 
all  second  maximal  subgroups  are  normal.
In  the  same  article  Huppert  proved  that  if  all  $3$-maximal 
subgroups  of  $G$  are  normal  in  $G,$
then  the  commutator  subgroup  $G'$
of  $G$  is    nilpotent   and  the 
chief   rank  of  $G$  is  at  most  $2.$
These  results  were  developed  by  many  authors. 
In particular, L.Ja. Poljakov \cite{polyakov} proved that $G$ is supersoluble 
if every $2$-maximal subgroup of $G$ is permutable
with every maximal subgroup of $G.$ He also established  the solubility of $G$
in the case when every maximal subgroup of $G$ permutes with every $3$-maximal 
subgroup of $G.$
Some later,  R.K. Agrawal \cite{agraval} proved that $G$ is supersoluble if
any  $2$-maximal subgroup of $G$ is permutable with every Sylow subgroup 
of $G.$
In \cite{janko1}, Z. Janko described the groups whose $4$-maximal 
subgroups are normal. A description of 
nonsoluble groups with all $2$-maximal subgroups nilpotent was obtained by  
M. Suzuki \cite{suzuki} and Z. Janko \cite{janko2}. 
In \cite{gagenjanko}, T.M. Gagen and Z. Janko gave a description of simple groups
whose $3$-maximal subgroups are nilpotent.
V.A. Belonogov \cite{belonogov} studied those groups in which every $2$-maximal 
subgroup is nilpotent. Continuing this, V.N. Semenchuk \cite{Sem}  obtained a 
description of soluble groups whose all $2$-maximal subgroups are 
supersoluble. 
A. Mann \cite{mann} studied the structure of the groups whose $n$-maximal subgroups 
are subnormal.
He proved that if all $n$-maximal subgroups of a soluble group  $G$ 
are subnormal and 
$|\pi (G)| \geq n + 1,$
then $G$ is nilpotent; but if  $|\pi (G)| \geq n - 1,$ 
then $G$ is  $\phi$-dispersive for some ordering  $\phi$ of 
the set of all primes.  Finally, in the case  $|\pi (G)| = n,$ 
Mann described $G$ completely. 
A.E. Spencer \cite{spencer} studied groups in which every $n$-maximal chain
contains a proper subnormal subgroup. 
In particular, Spencer proved that $G$ is a Schmidt group with abelian Sylow subgroups if
every $2$-maximal chain  of $G$
contains a proper subnormal subgroup.

Among  the 
recent  results  on  $n$-maximal  subgroups  
we  can  mention the paper of X.Y. Guo and K.P. Shum \cite{guo1}.
In this paper the authors proved that $G$ is soluble if all its $2$-maximal  subgroups  
enjoy  the  cover-avoidance  property.
W. Guo, K.P. Shum, A.N. Skiba and
B. Li   \cite{guo2,li,guo3} gave  
new  characterizations  of  supersoluble  groups  in  terms  of  
$2$-maximal  subgroups. 
Sh. Li  \cite{li1} obtained a 
classification  of  nonnilpotent  groups  whose  all  $2$-maximal  subgroups  are 
$TI$-subgroups.   
In \cite{legchek},  W. Guo, H.V. Legchekova and A.N. Skiba   described the groups  whose  every $3$-maximal 
subgroup permutes  with all maximal 
subgroups.   
In  \cite{lucenko1}, W. Guo, Yu.V. Lutsenko and A.N. Skiba gave a description of 
nonnilpotent  groups in  which  every  two  $3$-maximal  subgroups  are 
permutable.  
Yu.V. Lutsenko and A.N. Skiba \cite{lucenko2}
obtained a description of  the groups  whose  all  $3$-maximal  subgroups  are  
$S$-quasinormal.
Subsequently,  this  result  was  strengthened by Yu.V. Lutsenko and A.N. Skiba in  \cite{lucenko3}  to  provide 
a  description of the groups whose all $3$-maximal subgroups are subnormal. 
Developing some of the above-mentioned results, 
W. Guo, D.P. Andreeva and A.N.  Skiba  \cite{ags}
obtained a description of the groups in which 
every $3$-maximal chain contains a proper subnormal subgroup.
In  \cite{BES}, A. Ballester-Bolinches, L.M. Ezquerro and A.N. Skiba  obtained  
a full classification of the groups in which the second maximal subgroups of the Sylow 
subgroups cover or avoid the chief factors of some of 
its chief series.
In \cite{KM1}, V.N. Kniahina and V.S. Monakhov
studied those groups $G$ in which every $n$-maximal subgroup permutes with 
each Schmidt subgroup. In partiqular, it was be proved that if $n=1,2,3,$ 
then $G$ is metanilpotent; but if $n\geq 4$ and $G$ is soluble, then the 
nilpotent length of $G$ is at most 
$n-1.$
In \cite{KS1}, V.A. 
Kovaleva and A.N. Skiba
described the groups whose all $n$-maximal subgroups
are $\mathfrak U$-subnormal. In \cite{KS2},  the authors obtained a 
description 
of the groups with all $n$-maximal subgroups $\mathfrak F$-subnormal  for some 
saturated formation $\mathfrak F$.
In \cite{KM2}, V.S. 
Monakhov and V.N. Kniahina studied the groups 
with all $2$-maximal subgroups $\mathbb P$-subnormal.

Recall that  a subgroup $H$ of $G$ is said to be: 
(i) {\sl $\mathfrak U$-subnormal} in $G$ if there exists a chain of subgroups 
$H=H_0\leq H_1\leq \cdots \leq H_{n}=G$ 
such that 
$H_{i}/(H_{i-1})_{H_i}\in \mathfrak U,$ for $i=1,\ldots
,n;$
(ii) {\sl ${\mathfrak U}$-subnormal in the sense 
of Kegel} \cite{Keg} or  {\sl $K$-${\mathfrak U}$-subnormal} (see 
p. 236 in \cite{BE}) in $G$ if  there exists a chain of subgroups 
$H=H_{0}\leq  H_{1} \leq \cdots \leq H_{t}= G$ 
such that either $H_{i-1}$ is normal in $H_{i}$
or $H_{i}/(H_{i-1})_{H_{i}}\in {\mathfrak U}$  for all $i=1,  \ldots ,
t.$ 
It is evident that every subnormal subgroup is $K$-$\mathfrak 
U$-subnormal.
The inverse, in general, it is not true. For example, in the group $S_3$
a subgroup of order $2$ is $K$-$\mathfrak U$-subnormal and at the same time
it is not subnormal. This elementary observation and the results in 
\cite{lucenko3, ags, KS1, KS2} make natural the 
following questions: 

I. {\sl What is the structure of  $G$ under the condition that every
$2$-maximal subgroup of $G$ is $K$-$\mathfrak U$-subnormal?}

II. {\sl What is the structure of $G$ under the condition that every
$3$-maximal subgroup of $G$ is $K$-$\mathfrak U$-subnormal?}

Before continuing, recall that $G$ is called {\sl a minimal nonsupersoluble group}    provided
 $G$ does not belong to $\mathfrak U$ but every proper subgroup of $G$ belongs to
$\mathfrak U.$ 
Such groups were described by B. Huppert \cite{huppert} and K. Doerk \cite{D}.
We say that $G$ is a {\sl special Doerk-Huppert group} or an {\sl $SDH$-group} if $G$ is a minimal nonsupersoluble 
group such that
$G^{\mathfrak U}$ is a  minimal normal 
subgroup of $G.$

The solution of the first of the above-mentioned questions originates to \cite{KS1,KS2},
where, in particular, the following theorem was proved.

{\bf Theorem A$^*$.} {\sl Every $2$-maximal subgroup of $G$ is $\mathfrak U$-subnormal 
in $G$ if and only if $G$ is either supersoluble or an
$SDH$-group.}

If every $2$-maximal subgroup of $G$ is $K$-$\mathfrak U$-subnormal, then 
every maximal subgroup of $G$ is supersoluble (see Lemma 2.3 below).
Therefore in this case  $G$ is either supersoluble or a minimal nonsupersoluble group, 
hence $G$ is soluble by \cite{huppert}. Thus 
we get the following

{\bf Theorem A.} {\sl Every $2$-maximal subgroup of $G$ is $K$-$\mathfrak U$-subnormal 
in $G$ if and only if $G$ is either supersoluble or an
$SDH$-group.}

In this paper, on the bases of Theorem A 
we analyse  Question II. Note that since each subgroup of every supersoluble group is 
$K$-$\mathfrak U$-subnormal, we need, in fact, only consider the case when
$G$ is not supersoluble. 
But in this case, in view of \cite[Theorem A]{KS1} or \cite[Theorem A]{KS2},
$|\pi (G)|\leq 4.$ 
The following theorems are proved.

{\bf Theorem B} (See Theorem B in \cite{KS4}). {\sl 
Let $G$ be a nonsupersoluble group with $|\pi (G)|=2.$ Let $p,q$ 
be distinct prime divisors  
of $|G|$ and $G_p,$ $G_q$  be Sylow $p$-subgroup and $q$-subgroup 
of $G$ respectively.
Every $3$-maximal subgroup of $G$ 
is $K$-$\mathfrak U$-subnormal in $G$ if and only if
$G$ is a soluble group of one of the following 
types:}

I. {\sl $G$ is either a minimal nonsupersoluble group such that $|\Phi (G^{\mathfrak U})|$ is a prime or
an $SDH$-group.}

II. {\sl $G=G_p\rtimes G_q,$ where $G_p$ is the unique minimal normal subgroup of 
$G$ and every $2$-maximal subgroup of $G_q$ is an abelian group of 
exponent dividing  $p-1.$
Moreover, every maximal subgroup of $G$ containing $G_p$ is either supersoluble or 
an $SDH$-group and at least one of the maximal subgroups of $G$ is not supersoluble. }

III. {\sl $G=(G_p\times Q_1)\rtimes Q_2,$ where 
$G_q=Q_1\rtimes Q_2,$ 
$G_p$ and $Q_1$ are minimal normal subgroups of $G,$
$|Q_1|=q,$ $G_p\rtimes Q_2$ is an $SDH$-group and every maximal subgroup of 
$G$ containing
$G_p\rtimes Q_1$ is supersoluble. Moreover, if $p<q,$ then every $2$-maximal subgroup of $G$ is nilpotent.}

IV. {\sl $G=G_p\rtimes G_q,$ where $G_p$ is a minimal normal subgroup of 
$G$, $O_q(G)\ne 1,$ $\Phi (G)\ne 1,$ every maximal subgroup of $G$ containing $G_p$ is either supersoluble or 
an $SDH$-group and $G/\Phi (G)$ is a group of one of Types II or III.}

V. {\sl $G=(P_1\times P_2)\rtimes G_q,$ where  $G_p=P_1\times P_2,$ $P_1,$ 
$P_2$ are minimal normal subgroups of $G,$
every maximal subgroup of  $G$ containing $G_p$ is supersoluble,
$P_1\rtimes G_q$ is an  $SDH$-group and $P_2\rtimes G_q$ is either an  $SDH$-group or a supersoluble group with $|P_2|=p.$ }

VI. {\sl $G=G_p\rtimes G_q,$ where $\Phi (G_p)$ is a minimal normal subgroup of 
$G,$
every maximal subgroup of  $G$ containing $G_p$ is supersoluble and 
$\Phi (G_p)\rtimes G_q$ is an $SDH$-group.}

VII. {\sl Each of the subgroups $G_p$ and  $G_q$ is not normal in $G$  
and the following hold:}

(i) {\sl if $p<q,$ then $G=P_1\rtimes (G_q\rtimes P_2),$ where 
$G_p=P_1\rtimes P_2,$ $P_1$ is a minimal normal subgroup of 
$G,$ $|P_2|=p,$ $G_q=\langle a\rangle$ is a cyclic group and $\langle 
a^q\rangle$ is normal in $G.$
Moreover, $G$ has precisely three classes
of maximal subgroups whose representatives are  
$P_1\rtimes G_q,$ $G_q\rtimes P_2$ and $\langle a^q\rangle \rtimes G_p,$ 
where $P_1\rtimes G_q$ is an $SDH$-group;}

(ii) {\sl if $p>q,$ then $G=P_1(G_q\rtimes P_2),$ where $G_p=P_1P_2,$ $P_1$ is a normal subgroup of 
$G,$ $P_2=\langle b\rangle$ is a cyclic group and $1\ne P_1\cap P_2=\langle b^p\rangle .$ Moreover, $G$ 
has precisely three classes
of maximal subgroups whose representatives are 
$P_1\rtimes G_q,$ $G_q\rtimes P_2,$ $G_p,$ where
$|G:G_q\rtimes P_2|=p,$
$P_1\rtimes G_q$ is a normal supersoluble group and $G_q\rtimes P_2$ is an $SDH$-group.}

{\bf Theorem C.} {\sl 
Let $G$ be a nonsupersoluble group with $|\pi (G)|=3.$ Let $p,q,r$ 
be distinct prime divisors  
of $|G|$ and $G_p,$ $G_q$, $G_r$  be Sylow $p$-subgroup, $q$-subgroup 
and $r$-subgroup
of $G$ respectively.
Every $3$-maximal subgroup of $G$ 
is $K$-$\mathfrak U$-subnormal in $G$ if and only if
$G$ is a soluble group of one of the following 
types:}

I. {\sl $G$ is either a minimal nonsupersoluble group such that $|\Phi (G^{\mathfrak U})|$ is a prime or
an $SDH$-group.}

II. {\sl  $G=G_p\rtimes (G_q\rtimes G_r)$,  where $G_p$ is a minimal 
normal subgroup of $G$, every maximal subgroup of $G$ is either supersoluble or 
an $SDH$-group and at least one of the maximal subgroups of $G$ is not 
supersoluble. Moreover, the following hold:}

(i) {\sl if $G_p$ is the unique minimal normal subgroup of
$G$, then  every $2$-maximal subgroup of $G_q\rtimes G_r$ is an abelian group of 
exponent dividing  $p-1$;}
 
(ii) {\sl if $G_q\rtimes G_r$ is an $SDH$-group, then
all maximal subgroups of $G$ containing $G_pG_q$ are supersoluble and
$G_p\rtimes G_r$ is either an $SDH$-group or a supersoluble group with $|G_p|=p$.}

III. {\sl $G=(P_1\times P_2)\rtimes (G_q\rtimes G_r)$, where $G_p=P_1\times P_2$, $P_1$ and $P_2$ 
are minimal normal subgroups of $G$, $G_q$ and $G_r$ are cyclic groups, 
every maximal subgroup of  $G$ containing $G_p$ is supersoluble,
$P_1\rtimes (G_q\rtimes G_r)$ is an $SDH$-group and $P_2\rtimes (G_q\rtimes G_r)$ 
is either an $SDH$-group or a supersoluble group with $|P_2|=p$.}

IV. {\sl $G=G_p\rtimes (G_q\rtimes G_r)$, where $\Phi (G_p)$ is a minimal normal subgroup of $G$,
every maximal subgroup of  $G$ containing $G_p$ is supersoluble and 
$\Phi (G_p)\rtimes (G_q\rtimes G_r)$ is an $SDH$-group.}

{\bf Theorem D.} {\sl 
Let $G$ be a nonsupersoluble group with $|\pi (G)|=4.$ Let $p,q,r,t$ 
be distinct prime divisors  
of $|G|$ ($p>q>r>t$) and $G_p,$ $G_q$, $G_r$, $G_t$  be Sylow $p$-subgroup, $q$-subgroup, 
$r$-subgroup  and $t$-subgroup
of $G$ respectively.
Every $3$-maximal subgroup of $G$ 
is $K$-$\mathfrak U$-subnormal in $G$ if and only if
$G=G_p\rtimes (G_q\rtimes (G_r\rtimes G_t))$ is a soluble group such that $G$ 
has precisely three classes
of maximal subgroups whose representatives are  $G_qG_rG_t$, $G_pG_qG_r\Phi (G_t)$, 
$G_pG_q\Phi (G_r)G_t$, $G_p\Phi (G_q)G_rG_t$,  and
every nonsupersoluble maximal subgroup of $G$
is an $SDH$-group, $G_r$ and $G_t$ are cyclic groups and the following hold: }

(1) {\sl if $G_qG_rG_t$ is an $SDH$-group, then $G^{\mathfrak U}=G_p\times G_q$, $G_q$ is a minimal normal subgroup of $G$,
$G_pG_qG_r\Phi (G_t)$ and 
$G_pG_q\Phi (G_r)G_t$ are supersoluble and $G_pG_rG_t$ is either an $SDH$-group or a supersoluble group with $|G_p|=p$;}

(2) {\sl if $G_qG_rG_t$ is supersoluble, then $G_q$ is cyclic.}

Theorems B, C and D show that the class of the groups
with all $3$-maximal subgroups $K$-$\mathfrak U$-subnormal is essentially
wider then the class of the groups with all $3$-maximal subgroups 
subnormal \cite{lucenko3}.

All unexplained notation and terminology are standard. The reader is
referred to  \cite{BE}, \cite{DH} and \cite{Guo} if necessary.

\section{Preliminary Results}

Let $M$ be a maximal subgroup of $G$.
Recall that $M$ is said to be 
{\sl $\mathfrak U$-normal} in $G$ if $G/M_G\in \mathfrak U$,
otherwise it is said to be {\sl $\mathfrak U$-abnormal} in $G$.
Note that if $G$ is soluble, then
$M$ is $\mathfrak U$-normal in $G$ if and only if $|G:M|$ is a prime.

We use the following results.

{\bf Lemma 2.1.} 
{\sl Let  $H$ and $K$ be subgroups of $G$ such that $H$ is  $K$-$\mathfrak U$-subnormal in $G$. }

(1) {\sl $H\cap K$
is  $K$-$\mathfrak U$-subnormal in
$K$} \cite[6.1.7(2)]{BE}.

(2) {\sl If $N$ is a normal subgroup in $G$, then $HN/N$ is  
$K$-$\mathfrak U$-subnormal in $G/N$ }
\cite[6.1.6(3)]{BE}.

(3) {\sl If  $K$ is $K$-$\mathfrak U$-subnormal in $H$, then
$K$ is $K$-$\mathfrak U$-subnormal in $G$} \cite[6.1.6(1)]{BE}.

(4) {\sl If $G^{\mathfrak U}\leq K$, then $K$
is $K$-$\mathfrak U$-subnormal in $G$} \cite[6.1.7(1)]{BE}. 

The following lemma is evident.

{\bf Lemma 2.2.} {\sl If $G$ is supersoluble, then every subgroup of $G$ 
is $K$-$\mathfrak U$-subnormal in $G$. }

{\bf Lemma 2.3.} {\sl If every $n$-maximal subgroup of  $G$ is 
$K$-$\mathfrak U$-subnormal in $G,$ then every $(n-1)$-maximal subgroup of $G$
is supersoluble and 
every $(n+1)$-maximal subgroup of 
$G$ is $K$-$\mathfrak U$-subnormal in $G.$ }

{\bf Proof.} We first show that every $(n-1)$-maximal subgroup of $G$
is supersoluble. Let $H$ be an $(n-1)$-maximal subgroup of $G$
and  $K$ any maximal subgroup of $H.$
Then $K$ is an $n$-maximal subgroup of $G$ and so, by hypothesis,
$K$ is $K$-$\mathfrak U$-subnormal in $G.$ Hence $K$ 
is $K$-$\mathfrak U$-subnormal in $H$ by Lemma 2.1(1). 
Therefore either $K$ is normal in $H$ or  $H/K_H\in \mathfrak U.$ 
If $K$ is normal in $H,$ then $|H:K|$ is a prime in view of maximality of 
$K$ in $H.$ Let $H/K_H\in \mathfrak U.$ Then we also get that
$|H:K|=|H/K_H:K/K_H|$
is a prime.
Since $K$ is an arbitrary maximal subgroup of $H$, it follows that $H$ is supersoluble.

Now, let $E$ be an $(n+1)$-maximal subgroup of $G,$
and let $E_1$ and $E_2$ be an $n$-maximal and an
$(n-1)$-maximal subgroup of $G,$ respectively, such that
$E\leq E_1\leq E_2.$ Then, by the above, $E_2$ is supersoluble, so
$E_1$ is supersoluble. Hence  $E$ is 
$K$-$\mathfrak U$-subnormal in $E_1$ by Lemma 2.2. 
By hypothesis, $E_1$ is
$K$-$\mathfrak U$-subnormal in $G.$ Therefore $E$ is $K$-$\mathfrak U$-subnormal in $G$ by 
Lemma 2.1(3). The lemma is proved.

Fix some ordering $\phi$ of the set of all primes. 
The record $p\phi q$  means that $p$  precedes $q$ in $\phi$ and $p\ne q$.
Recall that 
a group $G$ of order 
$p_1^{\alpha _1}p_2^{\alpha _2}\cdots
p_n^{\alpha _n}$
is called {\sl $\phi$-dispersive} whenever 
$p_{1}\phi p_{2}\phi \cdots \phi p_{n}$
and for every $i$ there is 
a normal subgroup of $G$ of order 
$p_1^{\alpha _1}p_2^{\alpha _2}\cdots
p_i^{\alpha _i}$. 
Furthermore, if $\phi$ is such that $p\phi q$  always implies $p>q$, 
then every $\phi$-dispersive group is called  {\sl Ore dispersive}.

{\bf Lemma 2.4} (See \cite[Theorems ' and C]{KS1} or \cite[Theorems ' and D]{KS2}). {\sl Let $G$ be a soluble group 
in which every $n$-maximal subgroup is  
$K$-$\mathfrak U$-subnormal. If $|\pi (G)|\geq n$, then $G$ is  
$\phi$-dispersive for some ordering $\phi$  of the set of all primes. 
Moreover, if  $|\pi (G)|\geq n+1$, then $G$  
is Ore dispersive and the followig holds:
if $G$ is not supersoluble, then 
$G=A\rtimes B$, where $A=G^{\mathfrak U}$  and   $B$ are Hall subgroups of $G$, 
$A$ is either of the form $N_1\times \cdots \times N_t$, 
where each $N_i$ is a minimal normal subgroup of $G$, which is a Sylow 
subgroup of $G$,  for $i=1, \ldots , t$, or a Sylow $p$-subgroup of $G$ of exponent $p$
for some prime $p$.}

{\bf Lemma 2.5.} 
{\sl Let $G$ be a minimal nonsupersoluble group. The following hold:}

(1) {\sl $G$ is soluble and $|\pi (G)|\leq 3$} \cite{huppert};

(2) {\sl if $G$ is not a Schmidt group, then $G$ is Ore dispersive} 
\cite{huppert};

(3) {\sl $G^{\mathfrak U}$  is the unique normal Sylow subgroup of $G$} 
\cite{huppert, D};

(4) {\sl $G^{\mathfrak U}/\Phi (G^{\mathfrak U})$  is a noncyclic chief 
factor of $G$} \cite{D};

(5) {\sl if $S$ is a complement of $G^{\mathfrak U}$ in $G$, then
$S/S\cap \Phi (G)$ is either a primary cyclic group or a Miller-Moreno 
group} \cite{D}.

{\bf Lemma 2.6} (See \cite[2]{Sem}). {\sl Let $G$ be 
a minimal nonsupersoluble group with
$|\pi (G)|=3$. Let
$p_1$, $p_2$, $p_3$ be distinct prime divisors  
of $|G|$ such that $p_1>p_2>p_3$ and $G_{p_i}$ be a Sylow $p_i$-subgroup of 
$G$, $i=1,2,3$. The following hold:}

1) {\sl $G$ is Ore dispersive and $G_{p_2},G_{p_3}$ are cyclic;}

2) {\sl $G$ has precisely three classes
of maximal subgroups;}

3) {\sl $G_{p_1}/\Phi (G_{p_1}),G_{p_1}G_{p_2}/G_{p_1}\Phi 
(G_{p_2}),G_{p_1}G_{p_2}G_{p_3}/G_{p_1}G_{p_2}\Phi (G_{p_3})$ 
are chief factors of $G$.}

{\bf Lemma 2.7.} {\sl If every $3$-maximal subgroup of $G$ 
is $K$-$\mathfrak U$-subnormal in $G,$ then $G$ is soluble. }

{\bf Proof.} Suppose that lemma is false and let $G$ be a
counterexample with $|G|$ minimal. 
Since every $3$-maximal subgroup of $G$
is $K$-$\mathfrak U$-subnormal in $G,$ every $2$-maximal subgroup of $G$
is supersoluble by Lemma 2.3. Hence every maximal subgroup of  $G$ 
is either supersoluble or a minimal nonsupersoluble group. 
Therefore all proper subgroups of $G$ are soluble in view of Lemma 2.4(1).
Assume that all $3$-maximal subgroups of $G$ are identity. Then all  $2$-maximal subgroups of $G$ 
have prime orderes and so every maximal subgroup of $G$ 
is supersoluble. Hence $G$ is either supersoluble or a minimal nonsupersoluble group.
Thus in view of Lemma 2.4(1), $G$ is soluble, a contradiction.
Hence there is a $3$-maximal subgroup $T$ of $G$ such that $T\ne 1.$
Since $T$ is 
$K$-$\mathfrak U$-subnormal in $G,$ there exists a proper subgroup 
$H$ of $G$ such that $T\leq H$ and either $G/H_G\in \mathfrak U$ or $H$ is 
normal in $G.$
If $G/H_G\in \mathfrak U,$ then $G$ is soluble in view of solubility of $H_G,$ a contradiction.
Therefore $H$ is normal in  $G$.
Let $E/H$  be any  $3$-maximal subgroup of $G/H.$
Then  $E$  is a $3$-maximal subgroup of $G,$ hence $E$ 
is $K$-$\mathfrak U$-subnormal in  $G.$  
Hence $E/H$ is  $K$-$\mathfrak U$-subnormal in  $G/H$ 
by Lemma 2.1(2).  Thus the hypothesis holds
 for  $G/H.$ 
Hence $G/H$ is soluble by the choice of $G.$ Therefore $G$ is soluble. 
This contradiction completes the proof of 
the lemma.

\section{Proofs of Theorems C and D}

{\bf Proof of Theorem C.} {\sl Necessity.} By Lemma 2.7, $G$ is soluble.
Hence every $K$-$\mathfrak U$-subnormal subgroup of $G$ is $\mathfrak U$-subnormal in $G$.
Let $W$ be a maximal subgroup of $G$. In view of hypothesis and Lemma 
2.1(1), every $2$-maximal subgroup of $W$ is  $K$-$\mathfrak U$-subnormal in  
$W$. Therefore, by Theorem A, $W$ is either supersoluble or an  
$SDH$-group.
In particular, all $2$-maximal subgroups of $G$ 
are supersoluble.

If all maximal subgroups of  $G$ 
are supersoluble, then $G$ is a minimal nonsupersoluble group.
In view of Lemma 2.5(3),
$G^{\mathfrak U}=G_t$ is a Sylow $t$-subgroup of $G$ for some prime divisor $t$ of $|G|$. 
Suppose that $|\Phi (G_t)|\geq t^2$. 
Let $M$ be a maximal subgroup of $G$ such that $G_t\nleq M$. Then
$G=G_tM$ and $M=(G_t\cap M)G_{t'}=\Phi (G_t)G_{t'}$ by Lemma 2.5(4), where $G_{t'}$  
is a Hall $t'$-subgroup of $G$. Since $M$ is supersoluble, there is a $2$-maximal subgroup $E$ of $M$ 
such that $|M:E|=t^2$.
Hence $M=\Phi (G_t)E$ and so $G=G_tE$. 
Since $E$ is $\mathfrak U$-subnormal in $G$,  
there exists a proper subgroup $H$ of $G$ such that $E\leq H$ and $G/H_G\in 
\mathfrak U$. Therefore $G_t\leq H$, hence $G=G_tE\leq H$.
This contradiction shows that $|\Phi (G_t)|\leq t$. 
Thus $G$ is a group of Type I.

Now consider the case when at least one of the maximal subgroups of  $G$
is not supersoluble. In view of Lemma 2.4, 
$G$  is  $\phi$-dispersive for some ordering $\phi$  of the set of all 
primes. We can assume that $G=G_p\rtimes (G_q\rtimes G_r)$.

First suppose that $G_p$ is a minimal normal subgroup of $G$.
Then  $M=G_qG_r$ is a maximal subgroup of $G$.
Note that if $G_p$ is the unique minimal normal subgroup of $G$, then every
$2$-maximal subgroup of $M$ is an abelian group of exponent dividing
$p-1$. Indeed, in this case $F(G)=G_p$.
Let $K$ be any $2$-maximal subgroup of $G_qG_r$.
Then $G_pK$ is a $2$-maximal subgroup of $G$ and so $G_pK$
is supersoluble.
Hence $G_pK/F(G_pK)$ is an abelian group of exponent dividing $p-1$ by \cite[1, 1.5 and Appendixes, 3.2]{wein}. 
Since 
$C_G(G_p)=C_G(F(G))\leq F(G)=G_p$,  it follows that $C_G(G_p)=G_p$. Therefore 
$O_q(G_pK)=1$, hence $F(G_pK)=G_p$. Thus $K\simeq G_pK/G_p=G_pK/F(G_pK)$ 
is an abelian group of exponent dividing $p-1$.

Suppose that 
$M$ is an $SDH$-group. 
Then $G_q=M^{\mathfrak U}$ is a minimal normal subgroup of 
$M$. 
We show that every maximal subgroup $V$ of $G$ containing $G_pG_q$ 
is supersoluble. Suppose that $V$ is an $SDH$-group. Then $|\pi 
(V)|=2$. Indeed, if $|\pi (V)|>2$, then $|\pi 
(V)|=3$ in view of Lemma  2.5(1). Hence $G_q$ is cyclic by Lemma 2.6(1) and so
$G_qG_r$ is supersoluble, a contradiction.
Therefore $V=G_pG_q$. Note that $\Phi (G_pG_q)\leq \Phi (G)$
in view of normality of $G_pG_q$ in 
$G$. Since $G_p$ is a minimal normal subgroup of $G$ and $G_q$ is a 
minimal normal subgroup of $M=G_qG_r$, we have $\Phi (G)\cap G_p=\Phi (G)\cap G_q=1$. 
Consequently, $\Phi (G)\leq G_r$ and so $\Phi (V)=\Phi 
(G_pG_q)=1$. By Lemma 2.5(5), $G_q$ is either a primary cyclic group or a Miller-Moreno 
group. Since $G_q$ is a minimal normal subgroup of $M$,  $G_q$ is abelian
in view of solubility of $M$. Hence $G_q$ is a cyclic group.  
This contradiction completes the proof of supersolubility of $V$.

Since $G_q$ is a minimal normal subgroup of $M=G_qG_r$, $G_r$ is a maximal 
subgroup of $M$. Thus $G_pG_r$ is a maximal subgroup of $G$.
Finally, show that if $G_pG_r$ is supersoluble, then $|G_p|=p$. 
Suppose that $G_pG_r$ is supersoluble, but $|G_p|\geq p^2$. 
In view of solubility of $G$, it follows that $M$ is $\mathfrak 
U$-abnormal in $G$. Since $G_pG_r$ is supersoluble,
$G_r$ is a $k$-maximal subgroup of $G_pG_r$ for some 
$k\geq 2$.
Hence $G_r$ is a $(k+1)$-maximal subgroup of $G$ and so $G_r$ is $\mathfrak U$-subnormal in $G$ 
by hypothesis and Lemma 2.3. Therefore there is a proper subgroup $H$ of 
$G$ such that $G_r\leq H$ and $G/H_G\in \mathfrak U$. Consequently, 
$G^{\mathfrak U}\leq H$. 
Moreover, in view of hereditary of $\mathfrak U$, $M^{\mathfrak 
U}\leq G^{\mathfrak U}$. Hence $M=G_qG_r=M^{\mathfrak U}G_r\leq 
G^{\mathfrak U}G_r\leq H$. Since $H$ is a proper subgroup of $G$ and $M$ is a maximal subgroup of $G$,
it follows that $M=H$ is $\mathfrak U$-normal in $G$. This contradiction 
shows that $|G_p|=p$. 
Thus $G$ is a group of Type II.

Now suppose that $G_p$ is not a minimal normal subgroup of 
 $G$.  
Let $W$ be a maximal subgroup of $G$ such that $G_p\leq W$. 
Then $G_p$ is not a minimal normal subgroup of $W$. Therefore 
$W$ is not an $SDH$-group. Hence $W$ is supersoluble.

Let $\Phi (G_p)=1$. By Maschke's 
Theorem, $G_p=P_1\times P_2$, 
where  $P_1$ is a minimal normal subgroup of $G$ and $P_2$ is a normal 
subgroup of $G$. Then $M=P_2G_qG_r$ is a maximal subgroup of $G$. 
We show that  $P_2$ is also a minimal normal subgroup of
$G$. If $M$ is an $SDH$-group, then $P_2=M^{\mathfrak U}$ is 
a minimal normal subgroup of $M$,
so $P_2$ is a minimal normal subgroup of $G$. Assume that
$M$ is supersoluble.
Then $G/P_1\simeq M$ is a supersoluble group. If $P_1G_qG_r$ is supersoluble, 
then $G/P_2\simeq P_1G_qG_r$ is supersoluble and hence $G$ is supersoluble, a 
contradiction. Thus $P_1G_qG_r$ is not a supersoluble group. But every 
$2$-maximal subgroup of $G$ is supersoluble. Hence $P_1G_qG_r$
is a maximal subgroup of $G$, so $P_2$ is a minimal normal subgroup of $G$. 

Since $G$ is not supersoluble, at least one of the subgroups $M=P_2G_qG_r$ or 
$L=P_1G_qG_r$ is not supersoluble. Let 
$L$ be an $SDH$-group. 
Then in view of Lemma 2.6(1), $G_q$ 
and  $G_r$ are cyclic. We show also that if $M$ 
is supersoluble, then $|P_2|=p$.
Suppose that 
$|P_2|\geq p^2$. Then $L$ is $\mathfrak U$-abnormal in $G$. 
Note also that  $G^{\mathfrak U}=P_1$ since $G/P_1\simeq M$ is supersoluble.
Since $|P_2|\geq p^2$ and $M$ is supersoluble, 
$G_qG_r$ is a $k$-maximal subgroup of $M$, 
where $k\geq 2$. 
Hence $G_qG_r$ is a $(k+1)$-maximal subgroup of $G$ and so 
$G_qG_r$ is $\mathfrak 
U$-subnormal in $G$ by hypothesis and Lemma 2.3. 
Let $H$ be a proper subgroup of $G$ such that 
$G_qG_r\leq H$ and $G/H_G\in \mathfrak U$. Then $P_1=G^{\mathfrak 
U}\leq H$, hence $L=P_1G_qG_r\leq H$. Since $H$ is a proper subgroup of 
$G$ and $L$ is a maximal subgroup of $G$, it follows that $L=H$ is 
$\mathfrak U$-normal subgroup of $G$. This contradiction 
shows that $|P_2|=p$. Thus $G$ is a group of Type III.

Let now 
$\Phi (G_p)\ne 1$. 
By the above, all maximal subgroups of $G$ containing 
$G_p$ are supersoluble.
Since $G$ is not a minimal nonsupersoluble group,
there is a maximal subgroup $U$ of $G$ such that 
$G_p\nleq U$ and $U$ is an $SDH$-group. 
Let $U_p$ be a Sylow $p$-subgroup of $U$. Then $U_p=G_p\cap U$ 
is normal in $U$ and so $U_p=U^{\mathfrak U}$ 
is a minimal normal subgroup of $U$.
But since $1\ne \Phi (G_p)\leq U_p$ and $\Phi (G_p)$ is normal in $U$, 
it follows that $\Phi (G_p)=U_p$ is a minimal normal 
subgroup of $U$.
Therefore $\Phi (G_p)$ is a minimal normal 
subgroup of $G$.
Thus  
$G$ is a group of Type IV.

{\sl Sufficiency.} Let $E$ be any $3$-maximal subgroup of $G$. Let 
$M$ be a maximal subgroup of  $G$ such that $E$ is a $2$-maximal subgroup of $M$. 
Since $G$ is a group of one of Types I--IV, every maximal subgroup of $G$ is either 
supersoluble or an $SDH$-group. In particular, in view of Lemma 2.2 and Theorem A,  
$E$ is $K$-$\mathfrak U$-subnormal in $M$. If 
$G^{\mathfrak U}\leq M$, then in this case $E$ is 
$K$-$\mathfrak U$-subnormal in $G$ by Lemma 2.1(3)(4).

Suppose that $G^{\mathfrak U}\nleq M$. Let $D=G^{\mathfrak U}E$. 
In view of Lemma 2.1(4), $D$ 
is $K$-$\mathfrak U$-subnormal in $G$. If $G$ is a group of Type
I, then $D$ is supersoluble and so $E$ is $K$-$\mathfrak 
U$-subnormal in  $D$ by Lemma 2.2. Hence 
$E$ is $K$-$\mathfrak 
U$-subnormal in $G$ by Lemma 2.1(3).

Let $G$ be a  group of Type II. First assume that $G_qG_r$ is supersoluble.
Then $G^{\mathfrak U}=G_p$ and so $M=G/G^{\mathfrak U}\simeq G_qG_r$. 
Since $E$ is a $2$-maximal subgroup of $M$, 
$D=G^{\mathfrak U}E=G_pE$ is a  $2$-maximal subgroup of $G$. 
Since every maximal subgroup of $G$ is either supersoluble or an 
$SDH$-group, every  $2$-maximal subgroup of $G$
is supersoluble. Therefore $D$ is supersoluble. As abovei t follows  that 
$E$ is  $K$-$\mathfrak U$-subnormal in  $G$.

Now assume that $G_qG_r$ is an $SDH$-group. 
In view of Lemma 2.5(4), $|G_q|\geq q^2$. Hence $G_pG_r$ is an 
$\mathfrak U$-abnormal subgroup of $G$ in view of solubility of $G$.
Suppose that $G_pG_r$ is supersoluble. 
Then 
$|G_p|=p$, hence $G_qG_r$ is $\mathfrak U$-normal in $G$.
Consequently, $G^{\mathfrak U}\leq G_qG_r$.
It is easy to see that every maximal subgroup $W$ of $G$ such that  
$G_pG_q\leq W$ is normal in $G$. Therefore $W$ is $\mathfrak U$-normal in $G$.
Hence $M=G_pG_r^x$ for some $x\in G$.
In view of supersolubility of $M$, $E$ is one of the subgroups $R_1$ or $G_pR_2$,
where $R_1$ is a maximal subgroup of $G_r^x$ and $R_2$ is a $2$-maximal subgroup of $G_r^x$. 
Thus $D=G^{\mathfrak U}E\leq G_pG_qR$, where $R$ is a maximal subgroup of $G_r^x$. 
Since every maximal subgroup of $G$ containing $G_pG_q$ is supersoluble, $D$ 
is supersoluble. Consequently,  arguing as 
above we get that $E$ is $K$-$\mathfrak U$-subnormal in $G$.

Finally, suppose that  $G_pG_r$
is an $SDH$-group. 
Then $|G_p|\geq p^2$ by Lemma 2.5(4) and so $G_qG_r$ is an  
$\mathfrak U$-abnormal subgroup of $G$.
By the above, every maximal subgroup of $G$ containing  
$G_pG_q$ is $\mathfrak U$-normal in $G$.
Therefore we can assume that  $M$ is one of the subgroups  
$G_pG_r$ or $G_qG_r$.  Let $M=G_pG_r$. Then $E$ is one of the subgroups 
$R$, $PR$ or $PR_1$, where  $P$ is a maximal subgroup of $G_p$, 
$R$ is a maximal subgroup of $G_r$ and  $R_1$ is a
$2$-maximal subgroup of $G_r$. In all these cases we get that there is a maximal subgroup $W$ of $G$
such that $G_pG_q\leq W$ and $E\leq W$. Since $W$ is supersoluble, as above it follows  that 
$E$ is 
$K$-$\mathfrak U$-subnormal in $G$.  Arguing as 
above we get also that in the case when $M=G_qG_r$, $E$ is $K$-$\mathfrak U$-subnormal in $G$.

Let $G$ be a group of one of Types III or IV.
Then  
$G^{\mathfrak U}\leq G_p$.  In view of supersolubility of every maximal subgroup of $G$ containing 
$G_p$,
$D\leq G_pE$ is supersoluble. Therefore as 
above we get that $E$ is $K$-$\mathfrak U$-subnormal in $G$.
The theorem is proved.

{\bf Proof of Theorem D.} {\sl Necessity.} By Lemma 2.7, $G$ is soluble.
Hence every $K$-$\mathfrak U$-subnormal subgroup of $G$ is $\mathfrak U$-subnormal in $G$.
As in the proof of Theorem C  we get that every  maximal subgroup of $G$ is either supersoluble or  
an $SDH$-group.
In particular, all $2$-maximal subgroups of $G$ 
are supersoluble.
If all maximal subgroups of  $G$ 
are supersoluble, then $G$ is a minimal nonsupersoluble group and so 
$|\pi (G)|=3$ by Lemma 2.5(1). This contradiction shows that 
at least one of the maximal subgroups of $G$ is an  $SDH$-group.

By Lemma 2.4, $G$ is Ore dispersive. 
Therefore  
$G=G_p\rtimes (G_q\rtimes (G_r\rtimes G_t))$.
We show that $G_p$ is a minimal normal subgroup of $G$.
Suppose that $G_p$ is not a minimal normal subgroup of $G$.
Let  $M$ be any maximal subgroup of  $G$ such that  
$G_p\nleq M$.
Since $G_p$ is not a minimal normal subgroup of $G$,
$M\cap G_p\ne 1$ and so $|\pi (M)|=4$. Hence $M$ is supersoluble in view of Lemma 2.5(1).
Let $L$ be any maximal subgroup of $G$ containing $G_p$. If $L$ 
is an $SDH$-group, then $G_p=L^{\mathfrak U}$ is a minimal normal subgroup 
of $L$. It follows that $G_p$ is a minimal normal subgroup of $G$, a contradiction. 
Consequently, $L$ is supersoluble. Thus all maximal subgroups of $G$ 
are supersoluble.
This contradiction shows that $G_p$ is a minimal normal subgroup of $G$. 
In particular, it follows that all maximal subgroups of $G$ containing no 
$G_p$ are pairwise conjugate in $G$.

Let $W=G_qG_rG_t$. First suppose that $W$ is an $SDH$-group.
Then $G_q=W^{\mathfrak U}$ is a minimal normal subgroup of 
$W$. 
Since 
$G/G_pG_q\simeq G_rG_t\in \mathfrak U$ and $G/G_p\simeq G_qG_rG_t=W$ is 
not supersoluble, $G^{\mathfrak U}=G_pG_q$  in view of Lemma  2.4.  Moreover,  
$G_q$ is a minimal normal subgroup of $G$ by Lemma 2.4. 
Note also that  every maximal subgroup of $G$ containing $G_pG_q$ is supersoluble in view of Lemma 2.5(3).

By Lemma  2.6, $G_r$ and $G_t$ are cyclic and $W$ 
has precisely three classes
of maximal subgroups whose representatives are
$G_rG_t$, $G_qG_r\Phi (G_t)$ and $G_q\Phi 
(G_r)G_t$.
Hence $G$ has precisely three classes
of maximal subgroups containing $G_p$ whose representatives are
$G_pG_rG_t$, $G_pG_qG_r\Phi (G_t)$ and $G_pG_q\Phi 
(G_r)G_t$.
Suppose that $G_pG_rG_t$ is supersoluble. 
We show that in this case  $|G_p|=p$.  Indeed, if  $|G_p|\geq 
p^2$, then $G_rG_t$ is a 
$k$-maximal subgroup of  $G_pG_rG_t$ ($k\geq 2$) in view of supersolubility 
of $G_pG_rG_t$. Hence $G_rG_t$ 
is a $(k+1)$-maximal subgroup of $G$. Thus $G_rG_t$ is  
$\mathfrak U$-subnormal in  $G$ by hypothesis  and Lemma 2.3. 
Therefore there is a proper subgroup $H$ of $G$ such that
$G_rG_t\leq H$ and $G/H_G\in \mathfrak U$. Then $G_pG_q=G^{\mathfrak 
U}\leq H$ and so $G=G_pG_qG_rG_t\leq H$. This contradiction shows that $|G_p|=p$. 

Now suppose that  $W=G_qG_rG_t$ is supersoluble.
In this case $G^{\mathfrak U}=G_p$.
Since $G$ is not a minimal nonsupersoluble group, 
there exists a maximal subgroup $M$ of $G$ such that $G_p\leq M$  and $M$ 
is an $SDH$-group. 
Since $G^{\mathfrak U}=G_p\leq M$, $M$ is $\mathfrak 
U$-normal in $G$.
Hence $|G:M|$ is a prime. Moreover, in view of solubility of $G$ and Lemma 2.5(1), $|\pi (M)|=3$. 
If $|G:M|=t$, then by the above  $|G_t|=t$. Furthermore,  
$G_q$ and $G_r$ are cyclic by Lemma 2.6(1).
Arguing as 
above we get that
in the cases  $|G:M|=q$ and $|G:M|=r$ the subgroups 
$G_q$, $G_r$ and  $G_t$ are cyclic.
Hence $W$ has precisely three classes
of maximal subgroups whose representatives are
$\Phi (G_q)G_rG_t$, $G_qG_r\Phi (G_t)$ and $G_q\Phi 
(G_r)G_t$.
Therefore $G$ has precisely three classes
of maximal subgroups containing $G_p$ whose representatives are
$G_p\Phi (G_q)G_rG_t$, $G_pG_qG_r\Phi (G_t)$ and $G_pG_q\Phi 
(G_r)G_t$. 

{\sl Sufficiency.}  
Let $E$ be any $3$-maximal subgroup of $G$. Let 
$M$ be a maximal subgroup of  $G$ such that $E$ is a $2$-maximal subgroup of $M$. 
Since  every maximal subgroup of $G$ is either 
supersoluble or an $SDH$-group, 
$E$ is $K$-$\mathfrak U$-subnormal in $M$ by  Lemma 2.2 and Theorem A.
Hence in the case when $G^{\mathfrak U}\leq M$,  $E$ is 
$K$-$\mathfrak U$-subnormal in $G$ by Lemma 2.1(3)(4).

Suppose that  $G^{\mathfrak U}\nleq M$. Let $D=G^{\mathfrak U}E$. 
By Lemma 2.1(4), $D$ 
is $K$-$\mathfrak U$-subnormal in $G$. 

Assume that $G_qG_rG_t$ is an $SDH$-group. In this case $G^{\mathfrak U}=G_pG_q$
and so we can assume that $M$ is one of the subgroups $G_qG_rG_t$ or $G_pG_rG_t$. 
If $M=G_qG_rG_t$, then $M$ is an $SDH$-group.
Hence $|M:E|$ is divisible by at least one of the numbers $r$ or $t$. Therefore 
there is a maximal subgroup $W$ of $G$ such that $D=G^{\mathfrak U}E=G_pG_qE\leq W$.
By hypothesis, $W$ is supersoluble and so $D$ is supersoluble.
Consequently, $E$ is 
$K$-$\mathfrak U$-subnormal in $D$ by Lemma 2.2.
Therefore $E$ is
$K$-$\mathfrak U$-subnormal in $G$ by Lemma 2.1(3).
Let $M=G_pG_rG_t$. If $M$ is an $SDH$-group, then as above we have that $E$
is    
$K$-$\mathfrak U$-subnormal in  $G$.
Suppose that $M$ is supersoluble. Then $|G_p|=p$ and so  
$|M:E|$  is divisible by at least one of the numbers $r$ or $t$.
Therefore arguing as 
above we get that
$E$ is
$K$-$\mathfrak U$-subnormal in $G$.

Let now $G_qG_rG_t$ is supersoluble. Then
$G^{\mathfrak U}=G_p$. Therefore we can asume that $M=G_qG_rG_t$. Since $E$ is a 
$2$-maximal subgroup of $M$,  $D=G^{\mathfrak U}E=G_pE$ is a
$2$-maximal subgroup of $G$.
Since every maximal subgroup of $G$ is either supersoluble or an 
$SDH$-group, every  $2$-maximal subgroup of $G$
is supersoluble. Therefore $D$ is supersoluble and so as above we get that 
$E$ is  $K$-$\mathfrak U$-subnormal in $G$.
The theorem is proved.

\

Note that the classes of groups which are described in Theorems B, C, and D are pairwise disjoint.  
It is easy to construct  examples to show that all these classes are not empty.


\begin{thebibliography}{s2}

\bibitem{redei}
L. R\'{e}dei, "Ein Satz uber die endlichen einfachen Gruppen", \emph{Acta Math.}, {\bf 84} (1950), 129--153.

\bibitem{huppert}                                
B. Huppert, "Normalteiler and maximale Untergruppen
endlicher Gruppen", \emph{Math. Z.}, {\bf 60} (1954), 409--434.

\bibitem{polyakov}
L.Ja. Poljakov, "Finite groups with permutable
subgroups", in \emph{Proc. Gomel Sem.: Finite groups}, Nauka i Tekhnika, Minsk, 1966, 75--88.


\bibitem{agraval}
R.K. Agrawal, "Generalized center and hypercenter of a finite group", \emph{Proc. 
Amer. Math. Soc.}, {\bf 54} (1976), 13--21.



\bibitem{janko1} Z. Janko, "Finite groups with invariant fourth maximal 
subgroups", \emph{Math. Zeitschr.}, {\bf82} (1963),  82--89. 

\bibitem{suzuki}
M. Suzuki, "The nonexistence of a certain type of simple groups of odd 
order", \emph{Proc. Amer. Math. Soc.}, {\bf 8}(4) (1957), 686--695. 

\bibitem{janko2}
Z. Janko, "Endliche Gruppen mit lauter nilpotent zweitmaximalen Untergruppen", \emph{Math. Z.}, {\bf 79} (1962), 422--424. 

\bibitem{gagenjanko}
T.M. Gagen, Z. Janko, "Finite simple groups with nilpotent third maximal 
subgroups", \emph{J. Austral. Math. Soc.}, {\bf 6}(4) (1966), 466--469.

\bibitem{belonogov}
V.A. Belonogov, "Finite soluble groups with nilpotent $2$-maximal 
subgroups", \emph{Math. Notes}, {\bf 3}(1) (1968), 15--21.


\bibitem{Sem} V.N. Semenchuk, "Soluble groups with supersoluble second 
maximal subgroup", \emph{Voprosy Algebry}, {\bf 1} (1985), 86--96.

\bibitem{mann}  A. Mann, "Finite groups whose $n$-maximal 
subgroups are subnormal", \emph{Trans. Amer. Math. Soc.}, {\bf 132} (1968), 395--409. 

\bibitem{spencer}
A.E Spencer, "Maximal nonnormal chains in finite groups", \emph{Pacific J. 
of Math.}, {\bf 27}(1) (1968),  167--173.



\bibitem{guo1}
X.Y. Guo, K.P. Shum,   
"Cover-avoidance  properties  and  the  structure  of  finite  groups",  
\emph{J.  Pure  Appl.  Algebra}, {\bf 181} (2003), 297--308. 

\bibitem{guo2} W. Guo, K.P. Shum, A.N. Skiba, "$X$-Semipermutable subgroups 
of finite groups", \emph{J. Algebra}, {\bf 315} (2007), 31--41. 
 

\bibitem{li}
B. Li, A.N. Skiba, "New characterizations of finite 
supersoluble groups", \emph{Sci. China Ser. A: Math.}, {\bf 50}(1) (2008), 827--841. 
 

\bibitem{guo3}
W. Guo, A.N. Skiba, "Finite groups with given $s$-embedded and 
$n$-embedded subgroups", \emph{J. Algebra}, {\bf 321} (2009), 2843--2860. 


\bibitem{li1}
Sh. Li, "Finite  non-nilpotent  groups  all  of  whose  second  
maximal  subgroups  are $TI$-groups", \emph{Math.  Proc.  of  the 
Royal Irish Academy}, {\bf 100A}(1) (2000), 65--71. 

\bibitem{legchek} 
W. Guo, E.V. Legchekova, A.N. Skiba, "Finite groups in which every $3$-maximal subgroup commutes 
with all maximal subgroups", \emph{Math. Notes}, {\bf 86}(3--4) (2009), 325--332.


\bibitem{lucenko1}                                
W. Guo, Yu.V. Lutsenko, A.N. Skiba, "On nonnilpotent groups with every two $3$-maximal subgroups
permutable", \emph{Siberian Math. J.}, {\bf 50}(6) (2009), 988--997. 


\bibitem{lucenko2}                                
Yu.V. Lutsenko, A.N. Skiba,  "Structure of finite groups with $S$-quasinormal third maximal 
subgroups", \emph{Ukrainian Math. J.}, {\bf 61}(12) (2009), 1915--1922.

\bibitem{lucenko3} Yu.V. Lutsenko, A.N. Skiba,  "Finite groups with subnormal second 
and third maximal subgroups", \emph{Math. Notes}, {\bf 91}(5) (2012), 680--688. 


\bibitem{ags}  W. Guo, D.P. Andreeva, A.N. Skiba, "Finite groups of Spencer 
hight $\leq 3$", \emph{Algebra Colloquium}, (in  Press). 

\bibitem{BES} A. Ballester-Bolinches, L.M. Ezquerro, A.N. 
Skiba, "On second maximal subgroups of Sylow subgroups of finite groups", \emph{J. 
Pure Appl. Algebra}, {\bf 215}(4) (2011), 705--714. 

\bibitem{KM1} V.N. Kniahina, V.S. Monakhov, "On the permutability of 
$n$-maximal subgroups with Schmidt subgroups", \emph{Trudy Inst. Mat. i Mekh. UrO 
RAN}, {\bf 18}(3) (2012), 125--130.

\bibitem{KS1} V.A. Kovaleva, A.N. Skiba, "Finite solvable groups with all $n$-maximal 
subgroups $\mathfrak U$-subnormal", \emph{Sib. Math. J.}, {\bf 54}(1) 
(2013), 65--73.

\bibitem{KS2} V.A. Kovaleva, A.N. Skiba, "Finite soluble groups with all $n$-maximal 
subgroups $\mathfrak F$-subnormal", \emph{J. Group Theory},  {\bf 17} (2014), 273--290.

\bibitem{KM2} V.S. Monakhov, V.N. Kniahina,
"Finite groups with $\mathbb P$-subnormal subgroups",
\emph{Ricerche di Matematica}, {\bf 62}(2) (2013), 307--322. 



\bibitem{Keg} O.H. Kegel,  "Zur Struktur mehrfach faktorisierbarer endlicher
Gruppen", \emph{Math. Z.}, {\bf 87} (1965), 409--434.

\bibitem{BE} A. Ballester-Bolinches, L.M. Ezquerro,
\emph{Classes of Finite Groups}, Springer-Verlag, Dordrecht, 2006.


\bibitem{D} K. Doerk, "Minimal nicht uberauflosbare, endliche Gruppen", 
\emph{Math. Z.}, {\bf 91} (1966),  198--205.

\bibitem{KS4} Xiaolan Yi, V.A. Kovaleva, "Finite biprimary  groups with all $3$-maximal subgroups
$K$-$\mathfrak U$-subnormal", \emph{Acta Math. Hung."} (submitted).


\bibitem{DH}  K. Doerk, T. Hawkes,  \emph{Finite Soluble Groups},
Walter de Gruyter, Berlin-New York, 1992.



\bibitem{Guo} W. Guo,  \emph{The Theory of Classes of  Groups}, Science
Press-Kluwer Academic Publishers, Beijin-New York-Dordrecht-Boston-London, 2000.

\bibitem{wein} M. Weinstein (ed.), etc,  \emph{Between Nilpotent and Solvable}, Polygonal
Publishing House, Passaic N.J., 1982.






\end{thebibliography}
 \end{document}